\documentclass[12pt]{article}
\usepackage{amsmath,amsthm,amssymb,amscd}
\usepackage{latexsym}
\newtheorem{thm}{Theorem}[section]
 \newtheorem{cor}[thm]{Corollary}
 \newtheorem{lem}[thm]{Lemma}
 \newtheorem{prop}[thm]{Proposition}
 \theoremstyle{definition}
 
 \theoremstyle{remark}

 \numberwithin{equation}{section}

\begin{document}

\title
{Hessian estimates for the conjugate heat equation coupled with the Ricci flow}

\author{ Hong Huang}
\date{}
\maketitle
\begin{abstract}
  In this short note we obtain some local and global upper bounds  for the Hessian of a positive solution to the conjugate heat equation coupled with the Ricci flow. 



\end{abstract}


\section {Introduction}

In the literature there are extensive researches on  derivative estimates for  solutions to the heat equation or conjugate heat equation coupled with the Ricci flow; see  for example  \cite{BCP}, \cite{Ca}, \cite{CH}, \cite{CTY}, \cite{C+}, \cite{EKNT}, \cite{H}, \cite{HZ}, \cite{Hu}, \cite{KZ}, \cite{Li}, \cite{N}, \cite{P}, and \cite{Z06}.

In this short note we  obtain some local and global upper bounds  for the Hessian of a positive solution to the conjugate heat equation coupled with the Ricci flow. I was inspired by Han-Zhang \cite{HZ}. To state our results we first  introduce some notations.
 Fix $T >0$. Let $(M, g(t))$, $t\in [0,T]$, be a (not necessarily complete) solution to Hamilton's Ricci flow
 \begin{equation*}
 \frac{\partial g(t)}{\partial t}=-2Ric (g(t))
 \end{equation*}
 on a manifold $M$ (without boundary) of dimension $n$.  For any points $x, y \in M$, let $d(x, y,t)$ be the distance between   $x$ and $y$  w.r.t. $g(t)$. For $(x_0,t_0)\in M\times (0,T]$, $r>0$, and $T'>0$, let
 \begin{equation*}
 Q_{r,T'}(x_0,t_0):=\{(x,t)\in M\times [0,T]  \hspace{1mm} | \hspace{1mm}  d(x,x_0,t)\leq r, \hspace{1mm} t_0-T'\leq t \leq t_0\}
\end{equation*}
 be a parabolic cube.

  Consider a smooth positive solution $u$ to the conjugate heat equation $(\frac{\partial}{\partial t}+\Delta_{g(t)}-R_{g(t)})u=0$ coupled with the Ricci flow on $M\times [0,T]$, where $R_{g(t)}$ is the scalar curvature of the metric $g(t)$. We have the following global upper bound for  $\text{ Hess} \hspace{0.5mm} u$.
\begin{thm} \label{thm 1.1} \    Let $(M^n,g(t))$, $t\in [0,T]$, be a Ricci flow on a compact manifold,  $u$ be  a positive solution to the conjugate heat equation coupled with the Ricci flow on $M\times [0,T]$ with $0 < u \leq A$. Then we have
\begin{equation*}
\nabla^2 u \leq u(\frac{18}{T-t}+C)(1+\log \frac{A}{u})g(t)   \hspace{8mm}  on  \hspace{2mm}  M\times [0,T),
\end{equation*}
and, in particular,
\begin{equation*}
\Delta u \leq u(\frac{18n}{T-t}+C)(1+\log \frac{A}{u})   \hspace{8mm}  on  \hspace{2mm}  M\times [0,T),
\end{equation*}
where $C$ depends  on $n$, the upper bounds of $|Rm|$, $|\nabla Ric|$, and $|\nabla^2 R|$ on $M\times [0,T]$.
\end{thm}

We also get a local upper bound for  $\text{ Hess} \hspace{0.5mm} u$.
\begin{thm} \label{thm 1.2}    Let $(M^n,g(t)$, $t\in [0,T]$,  be a (not necessarily complete) Ricci flow,  $u$ be  a positive solution to the conjugate heat equation coupled with the Ricci flow on $M\times [0,T]$ with $0 < u \leq A$. Let $x_0\in M$. Assume that  the parabolic cube  $Q_{4r,T}(x_0,T)$ is compact. Then we have
\begin{equation*}
\nabla^2 u \leq u(\frac{C_0}{T}+\frac{C_0}{r^2}+C_2)(1+\log \frac{A}{u})^2g(t)   \hspace{8mm}  on  \hspace{2mm}  Q_{r,\frac{T}{2}}(x_0,\frac{T}{2}),
\end{equation*}
where $C_0$ is a universal constant, and $C_2$ depends  on $n$, the upper bounds of $|Rm|$, $|\nabla Ric|$, and $|\nabla^2 R|$  on $Q_{4r,T}(x_0,T)$.
\end{thm}

In Section 2 we  derive local gradient  estimates  for a positive solution to the conjugate heat equation coupled with the Ricci flow.  In Section 3 we prove Theorems 1.1 and 1.2.

\vspace*{0.4cm}
\section{ Gradient estimates}

Let $(M^n,g(t))$, $t\in [0,T]$,  be a (not necessarily complete) Ricci flow, and $u$ be  a positive solution to the conjugate heat equation
\begin{equation*}
(\frac{\partial}{\partial t}+\Delta_{g(t)}-R_{g(t)})u=0
\end{equation*}
coupled with the Ricci flow on $M\times [0,T]$. Sometimes we'll  write $\partial_t$ for $\frac{\partial}{\partial t}$, and $u_t$ for $\frac{\partial u}{\partial t}$. Let $f=\log u$. Then
\begin{equation*}
f_t=-\Delta f-|\nabla f|^2+R.
\end{equation*}
Fix $\alpha >1$.
Let $\tau=T-t$, and
\begin{equation*}
F=\tau(\frac{|\nabla u|^2}{u^2}+\alpha \frac{u_t}{u}-\alpha R)=\tau(|\nabla f|^2+\alpha f_t-\alpha R).
\end{equation*}
Let $x_0\in M$. Assume that  the parabolic cube  $Q_{2r,T}(x_0,T)$ is compact, and  $-K_0\leq Ric \leq K_0$, $|\nabla R|\leq K_1$, and $\Delta R \leq K_2$ on $Q_{2r,T}(x_0,T)$.

\begin{lem} \label{lem 2.1}  With the above assumption, for any $\varepsilon >0$  we have
\begin{equation*}
\begin{split}
(\Delta+\partial_t) F & \geq -2\langle\nabla f, \nabla F\rangle +\frac{2\tau}{n+\varepsilon}(f_t+|\nabla f|^2-R)^2-(|\nabla f|^2+\alpha f_t-\alpha R)\\
& -2\tau |2-\alpha|K_0|\nabla f|^2-2\tau (\alpha-1)K_1|\nabla f|-\frac{n(n+\varepsilon)}{2\varepsilon}\alpha^2\tau K_0^2-\alpha \tau K_2 \\
\end{split}
\end{equation*}
on $Q_{2r,T}(x_0,T)$.
\end{lem}
\noindent {\bf Proof}. We have
\begin{equation*}
\begin{split}
\Delta |\nabla f|^2 &=2|\nabla^2f|^2+ 2\langle \nabla f, \nabla \Delta f\rangle+ 2Ric(\nabla f, \nabla f),\\
\Delta f_t &= (\Delta f)_t-2\langle Ric, \nabla^2f\rangle, \\
\tau \Delta f&=\tau(\alpha-1)(f_t-R)-F, \\
(|\nabla f|^2)_t &=2\langle \nabla f, \nabla f_t \rangle+2Ric(\nabla f, \nabla f \rangle, \\
\end{split}
\end{equation*}
and
\begin{equation*}
\begin{split}
& \Delta F \\
 =& \tau[2|\nabla^2f|^2+2\langle \nabla f, \nabla \Delta f\rangle+2 Ric(\nabla f, \nabla f)+\alpha ((\Delta f)_t-2\langle Ric, \nabla^2f\rangle)-\alpha \Delta R]\\
=& 2\tau |\nabla^2f|^2-2\langle\nabla f, \nabla F\rangle+2\tau(\alpha-1)\langle \nabla f, \nabla f_t \rangle-2\tau(\alpha-1)\langle \nabla R, \nabla f\rangle\\
& +2\tau Ric(\nabla f, \nabla f)+\alpha\tau(-f_{tt}-(|\nabla f|^2)_t+R_t-2\langle Ric, \nabla^2f\rangle)-\alpha\tau \Delta R \\
=& 2\tau |\nabla^2f|^2-2\langle\nabla f, \nabla F\rangle-F_t-\frac{F}{\tau}-2\alpha\tau \langle Ric, \nabla^2 f\rangle +2\tau(2-\alpha)Ric(\nabla f, \nabla f)\\
& -2\tau(\alpha-1)\langle \nabla R, \nabla f\rangle -\alpha \tau\Delta R. \\
\end{split}
\end{equation*}
See p. 61 of \cite{CTY}.

Note that
\begin{equation*}
|\langle Ric, \nabla^2f\rangle| \leq |Ric||\nabla^2f| \leq \frac{\alpha (n+\varepsilon)}{4\varepsilon}|Ric|^2 + \frac{\varepsilon}{\alpha (n+\varepsilon)}|\nabla^2f|^2
\end{equation*}
for any $\varepsilon>0 $, and
\begin{equation*}
|\nabla^2f|^2 \geq \frac{1}{n}(\Delta f)^2=\frac{1}{n}(f_t+|\nabla f|^2-R)^2,
\end{equation*}
 then the lemma follows. \hfill{$\Box$}

\begin{prop} \label{prop 2.2}  Let  $(M^n,g(t))$, $t\in [0,T]$, be a (not necessarily complete) Ricci flow and $u$ be a positive solution to the conjugate heat equation coupled with the Ricci flow on $M\times [0,T]$. Let $x_0\in M$. Suppose that   $Q_{2r,T}(x_0,T)$ is compact and $-K_0\leq Ric \leq K_0$, $|\nabla R|\leq K_1$, $\Delta R \leq K_2$ on  $Q_{2r,T}(x_0,T)$.
For any $\alpha>1$ and $\varepsilon>0$, we have
\begin{equation*}
\frac{|\nabla u|^2}{u^2}+\alpha \frac{u_t}{u}-\alpha R  \leq \frac{(n+\varepsilon)\alpha^2}{2(T-t)}+C(r^{-2}+r^{-1}+1)     \hspace{6mm}  on  \hspace{2mm}  Q_{r,T}(x_0,T)\setminus \{(x,T)  \hspace{1mm} |  \hspace{1mm} x\in M \},
\end{equation*}
where the constant $C$ depends  on $n,\alpha, \varepsilon, K_0,K_1$ and $K_2$.
\end{prop}
\noindent {\bf Proof}.
As in \cite{LY}, choose a smooth cutoff function $\psi: [0,\infty)\rightarrow [0,1]$ with $\psi=1$ on the interval $[0,1]$, $\psi=0$ on $[2,\infty)$, and
\begin{equation*}
 \psi'\leq 0,  \hspace{2mm} |\psi'|^2 \leq C_0  \psi, \hspace{2mm} \psi''\geq-C_0,
\end{equation*}
where $C_0$ is a universal constant.
Let $\phi(x,t)=\psi(\frac{d(x,x_0,t)}{r})$. Suppose that the maximum of the function $\phi F$ is positive, otherwise the result follows trivially. Assume that  $\phi F$ achieves its positive maximum at the point $(x_1, t_1)$. Then $(x_1,t_1) \in Q_{2r,T}(x_0,T)$ with $t_1 \neq T$. By  Calabi's trick \cite{C} we may assume that $\phi F$ is smooth at $(x_1, t_1)$. Let $\tau_1=T-t_1$. We compute at the point $(x_1, t_1)$ using Lemma 2.1,

\begin{equation*}
\begin{split}
0  \geq &  (\Delta+\partial_t) (\phi F) \\
\geq & \tau_1\phi \frac{2}{n+\varepsilon}(f_t+|\nabla f|^2-R)^2-CF\sqrt{\phi}r^{-1}|\nabla f|-\phi \frac{F}{\tau_1}\\
& -C\tau_1\phi |\nabla f|^2-C\tau_1 \phi-CF(r^{-2}+r^{-1}+1), \\
\end{split}
\end{equation*}
where the constant $C$ depends  on $n,\alpha, \varepsilon, K_0,K_1$ and $K_2$; compare \cite{CTY}, \cite{Li}, and \cite{S}.  Then we proceed as in \cite{LY} and \cite{CTY}.   \hfill{$\Box$}

\vspace*{0.4cm}

The following corollary is a slight improvement of Lemma 4.1 in \cite{CTY} in the Ricci flow case.
\begin{cor} \label{cor 2.3}  Suppose that $(M^n,g(t))$, $t\in [0,T]$, is a complete Ricci flow with $-K_0\leq Ric \leq K_0$, $|\nabla R|\leq K_1$, $\Delta R \leq K_2$ on $M\times [0,T]$, and that $u$ is a positive solution to the conjugate heat equation coupled with the Ricci flow on $M\times [0,T]$.
For any $\alpha>1$ and $\varepsilon>0$, at $(x,t)\in M\times [0,T)$  we have
\begin{equation*}
\frac{|\nabla u|^2}{u^2}+\alpha \frac{u_t}{u}-\alpha R \leq \frac{(n+\varepsilon)\alpha^2}{2(T-t)}+C,
\end{equation*}
where the constant $C$ depends  on $n,\alpha, \varepsilon, K_0,K_1$ and $K_2$.
\end{cor}

\begin{prop} \label{prop 2.4} Let  $(M^n,g(t))$, $t\in [0,T]$, be a (not necessarily complete) Ricci flow and $u$ be a positive solution to the conjugate heat equation coupled with the Ricci flow on $M\times [0,T]$. Let $x_0\in M$. Suppose that  $Q_{2r,T}(x_0,T)$ is compact and $-K_0\leq Ric \leq K_0$, $|\nabla R|\leq K_1$, $\Delta R \leq K_2$ on  $Q_{2r,T}(x_0,T)$.
For any $\alpha>1$ and $\varepsilon>0$, at $(x,t) \in Q_{r,T}(x_0,T)$ with $t\neq T$ we have
\begin{equation*}
\begin{split}
\frac{|\nabla u|^2}{u^2}&+\alpha \frac{u_t}{u}-\alpha R  \leq \frac{n\alpha^2}{T-t}+\frac{C\alpha^2}{r^2}(r\sqrt{K_0}+\frac{\alpha^2}{\alpha-1})+C\alpha^2K_0\\
& +\frac{n\alpha^2}{\alpha-1}(|2-\alpha| K_0+\frac{\alpha-1}{2}K_1)+n\alpha^2K_0+\alpha \sqrt{n(\alpha-1)K_1}+\alpha\sqrt{n\alpha K_2},
\end{split}
\end{equation*}
where the constant $C$ depends only on $n$.
\end{prop}
\noindent {\bf Proof}.  In Lemma 2.1 we let $\varepsilon=n$ and get
\begin{equation*}
\begin{split}
(\Delta+\partial_t) F & \geq -2\langle\nabla f, \nabla F\rangle +\frac{\tau}{n}(f_t+|\nabla f|^2-R)^2-(|\nabla f|^2+\alpha f_t-\alpha R)\\
& -\tau (2|2-\alpha|K_0+(\alpha-1)K_1)|\nabla f|^2-n\alpha^2\tau K_0^2-(\alpha-1)\tau K_1-\alpha \tau K_2. \\
\end{split}
\end{equation*}
Then we proceed as in \cite{Li} and \cite{S}.   \hfill{$\Box$}

\section{Hessian estimates}

Let $(M^n, g(t))$, $t\in [0,T]$, be a Ricci flow and $u$ be a  smooth positive solution to the conjugate heat equation  coupled with the Ricci flow on $M\times [0,T]$.

\begin{lem}\label{lem 3.1}  Let $u$ be  a positive solution to the conjugate heat equation coupled with the Ricci flow with $0 < u \leq A$, and $f = \log \frac{u}{A}$. Let  $u_{ij}$  denote  the Hessian of $u$ in local coordinates, and $v_{ij}:=\frac{u_{ij}}{u(1-f)}$.    In  local coordinates we have
\begin{equation*}
\begin{split}
(\partial_t+\Delta &-\frac{2f}{1-f}\nabla f\cdot \nabla)v_{ij} =\frac{|\nabla f|^2+Rf}{1-f}v_{ij}+\frac{1}{u(1-f)}[2R_{kijl}u_{kl}+R_{il}u_{jl}+R_{jl}u_{il} \\
& +2(\nabla_iR_{jl}+\nabla_jR_{il}-\nabla_lR_{ij})\nabla_lu -u\nabla_i\nabla_jR-\nabla_iR\nabla_ju-\nabla_jR\nabla_iu-Ru_{ij}]. \\
\end{split}
\end{equation*}
\end{lem}
\noindent Here we adopt the usual convention on summations.  For example, $R_{kijl}u_{kl}$ means $g^{ab}g^{pq}R_{aijp}u_{bq}$.  Note also that our convention on the curvature tensor $R_{ijkl}$ is the same as that of R. Hamilton, but is different from that of Han-Zhang \cite{HZ}.

\vspace*{0.4cm}

\noindent {\bf Proof}. The proof is similar to that of Lemma 3.1 in \cite{HZ}.   \hfill{$\Box$}

\begin{lem}\label{lem 3.2}  Let $u$ be  a positive solution to the conjugate heat equation coupled with the Ricci flow with $0 < u \leq A$, and $f = \log \frac{u}{A}$. For a function $h$ we let $h_i$ denote the 1-form $dh$ in local coordinates. We also let $w_{ij}:=\frac{u_iu_j}{u^2(1-f)^2}$ denote the 2-tensor $\frac{du \otimes du}{u^2(1-f)^2}$ in local coordinates.  In  local coordinates we have
\begin{equation*}
\begin{split}
(\partial_t+\Delta -\frac{2f}{1-f}\nabla f\cdot \nabla)w_{ij}& =\frac{2(|\nabla f|^2+Rf)}{1-f}w_{ij}+\frac{(Ru)_iu_j+u_i(Ru)_j}{u^2(1-f)^2} \\
& + 2(v_{ki}+fw_{ki})(v_{kj}+fw_{kj}) +R_{ik}w_{kj}+R_{jk}w_{ki}, \\
\end{split}
\end{equation*}
where $v_{ij}$ is as in Lemma 3.1.
\end{lem}
\noindent {\bf Proof}. The proof is similar to that of Lemma 3.2 in \cite{HZ}.   \hfill{$\Box$}

\vspace*{0.4cm}

\noindent {\bf Proof of Theorem 1.1}.
Let $V=(v_{ij})$, $W=(w_{ij})$, $w=\text{tr}\hspace{0.5mm} W=g^{ij}w_{ij}=\frac{|\nabla f|^2}{(1-f)^2}$,
and
\begin{equation*}
L=\partial_t+\Delta -\frac{2f}{1-f}\nabla f\cdot \nabla.
\end{equation*}
Then by Lemmas 3.1 and 3.2 we have
\begin{equation*}
\begin{split}
LV&=(1-f)wV+P, \\
LW&=2(1-f)wW+2(V+fW)^2+Q, \\
\end{split}
\end{equation*}
where $P$ is a 2-tensor whose $(i,j)$-th component in local coordinates is given by
\begin{equation*}
\begin{split}
P_{ij} & =\frac{Rf}{1-f}v_{ij}+\frac{1}{u(1-f)}[2R_{kijl}u_{kl}+R_{il}u_{jl}+R_{jl}u_{il} \\
& +2(\nabla_iR_{jl}+\nabla_jR_{il}-\nabla_lR_{ij})\nabla_lu -u\nabla_i\nabla_jR-\nabla_iR\nabla_ju-\nabla_jR\nabla_iu-Ru_{ij}], \\
\end{split}
\end{equation*}
and $Q$ is  a 2-tensor whose $(i,j)$-th component in local coordinates is given by
\begin{equation*}
Q_{ij} =\frac{2Rf}{1-f}w_{ij}+\frac{(Ru)_iu_j+u_i(Ru)_j}{u^2(1-f)^2} +R_{ik}w_{kj}+R_{jk}w_{ki}.
\end{equation*}

Now with the help of Theorem 10 in \cite{EKNT} and Corollary 2.3 here, we can proceed as in Han-Zhang \cite{HZ}.    \hfill{$\Box$}

\vspace*{0.4cm}

\noindent {\bf Proof of Theorem 1.2}.
 With the help of Theorem 10 in \cite{EKNT}, Proposition  2.2 here, and a space-time cutoff function, we can proceed as in Han-Zhang \cite{HZ} with minor modifications.   \hfill{$\Box$}

\begin{cor}\label{cor 3.3}  Let $(M^n,g(t)$, $t\in [0,T]$,  be a (not necessarily complete) Ricci flow,  $u$ be  a positive solution to the conjugate heat equation coupled with the Ricci flow on $M\times [0,T]$ with $0 < u \leq A$. Let $x_0\in M$. Assume that  the parabolic cube  $Q_{4r,T}(x_0,T)$ is compact. Then we have
\begin{equation*}
\nabla^2 u \leq u(\frac{C_0}{T-t}+\frac{C_0}{r^2}+C_2)(1+\log \frac{A}{u})^2g(t)   \hspace{8mm}  on  \hspace{2mm}  Q_{r,T}(x_0,T)\setminus \{(x,T)  \hspace{1mm} |  \hspace{1mm} x\in M \},
\end{equation*}
and, in particular,
\begin{equation*}
\Delta u \leq u(\frac{C_0n}{T-t}+\frac{C_0n}{r^2}+C_2)(1+\log \frac{A}{u})^2   \hspace{8mm}  on  \hspace{2mm}  Q_{r,T}(x_0,T)\setminus \{(x,T)  \hspace{1mm} |  \hspace{1mm} x\in M \},
\end{equation*}
where $C_0$ is a universal constant, and $C_2$ depends  on $n$, the upper bounds of $|Rm|$, $|\nabla Ric|$, and $|\nabla^2 R|$  on $Q_{4r,T}(x_0,T)$.
\end{cor}

\vspace*{0.4cm}

\noindent {\bf Acknowledgements}.   I would like to thank Professor Qi S. Zhang for helpful communications.  I was partially supported by NSFC 12271040 and Beijing Natural Science Foundation (Z190003).


\hspace *{0.4cm}

\vspace *{0.4cm}

Laboratory of Mathematics and Complex Systems (Ministry of Education),

School of Mathematical Sciences, Beijing Normal University,

Beijing 100875, P.R. China

 E-mail address: hhuang@bnu.edu.cn

\end{document}